\documentstyle[12pt]{article}
\newcommand{\sect}[1]{\section{#1}\setcounter{equation}{0}}

\font\mbn=msbm10 scaled \magstep1
\font\mbs=msbm7 scaled \magstep1
\font\mbss=msbm5 scaled \magstep1
\newfam\mbff
\textfont\mbff=\mbn
\scriptfont\mbff=\mbs
\scriptscriptfont\mbff=\mbss\def\mbf{\fam\mbff}
\def\Re{{\mbf R}}

\def\Z{{\mbf Z}}
\def\Co{{\mbf C}}

\def\Di{{\mbf D}}
\def\Bo{{\mbf B}}
\def\N{{\mbf N}}
\newtheorem{Th}{Theorem}[section]
\newtheorem{Lm}[Th]{Lemma}

\newtheorem{D}[Th]{Definition}
\newtheorem{Proposition}[Th]{Proposition}
\newtheorem{R}[Th]{Remark}

\author{Alexander Brudnyi\thanks{Research supported in part by NSERC and by
Max-Planck-Institut f\"{u}r Mathematik.\newline
}\\
Department of Mathematics and Statistics\\
University of Calgary, Calgary\\
Canada}
\title{Corona Theorem for $H^{\infty}$
on Coverings of Riemann Surfaces of
Finite Type}
\date{} 
\begin{document} 
\maketitle
\begin{abstract}
{In this paper continuing our work started in [Br1]-[Br3]
we prove the corona theorem for the algebra of bounded
holomorphic functions defined on an unbranched covering of a Caratheodory
hyperbolic Riemann surface of finite type.} 
\end{abstract}
\sect{\hspace*{-1em}. Introduction.}
\quad {\bf 1.1.}
Let $X$ be a complex manifold and let $H^{\infty}(X)$ be the Banach algebra 
of bounded holomorphic functions on $X$
equipped with the supremum norm. We assume that $X$ is Caratheodory
hyperbolic, that is, the functions in 
$H^{\infty}(X)$ separate the points of $X$. 
The maximal ideal space ${\cal M}={\cal M}(H^{\infty}(X))$ is the set of all
nonzero linear multiplicative functionals on $H^{\infty}(X)$. Since the norm
of each $\phi\in {\cal M}$ is $\leq 1$, ${\cal M}$
is a subset of the closed unit ball of the
dual space $(H^{\infty}(X))^{*}$. It is a
compact Hausdorff space in the Gelfand topology (i.e., in the weak $*$ 
topology induced by $(H^{\infty}(X))^{*}$). Further, there is a
continuous embedding $i:X\hookrightarrow {\cal M}$ taking $x\in X$ to the 
evaluation homomorphism $f\mapsto f(x)$, $f\in H^{\infty}(X)$. 
The complement to
the closure of $i(X)$ in ${\cal M}$ is called the {\em corona}.
The {\em corona problem} is: given $X$ to determine whether the corona is 
empty. For example, according to Carleson's celebrated Corona Theorem
[C] this is true for $X$ being the open unit disk in $\Co$. (This was
conjectured by Kakutani in 1941.) Also, there are
non-planar Riemann surfaces for which the corona is non-trivial (see, e.g.,
[JM], [G], [BD], [L] and references therein). The general problem for planar
domains is still open, as is the problem in several variables for the ball
and polydisk. (In fact, there are no known examples of domains in $\Co^{n}$,
$n\geq 2$, without corona.) At present, the strongest
corona theorem for planar domains is due to Jones and 
Garnett [GJ]. It states that the corona is empty for any 
Denjoy domain, i.e.,
a domain of the form $\overline{\Co}\setminus E$ where $E\subset\Re$. 

The corona problem has the following analytic reformulation, see, e.g., [Ga]:

A collection $f_{1},\dots, f_{n}$ of functions from $H^{\infty}(X)$ satisfies
the {\em corona condition} if
\begin{equation}\label{e1}
1\geq\max_{1\leq j\leq n}|f_{j}(x)|\geq\delta>0\ \ \ {\rm for\ all}\ \ \
x\in X.
\end{equation}
The corona problem being solvable (i.e., the corona is empty)
means that the Bezout equation
\begin{equation}\label{e2}
f_{1}g_{1}+\cdots+f_{n}g_{n}\equiv 1
\end{equation}
has a solution $g_{1},\dots, g_{n}\in H^{\infty}(X)$ for any 
$f_{1},\dots, f_{n}$ satisfying the corona condition. We refer to
$\max_{1\leq j\leq n}||g_{j}||_{\infty}$ as a ``bound on the corona
solutions``. (Here $||\cdot||_{\infty}$ is the norm on $H^{\infty}(X)$.)

The present paper is concerned with the corona problem for 
coverings of Riemann surfaces of finite type. Let us recall that a 
Riemann surface $Y$ is of finite type if the fundamental group $\pi_{1}(Y)$
is finitely generated. Our main result extends the class of Riemann surfaces
for which the corona theorem is true:
\begin{Th}\label{te1}
Let $r:X\to Y$ be an unbranched covering of a Caratheodory  hyperbolic
Riemann surface of finite type $Y$. Then $X$ is Caratheodory hyperbolic and
for any $f_{1},\dots, f_{n}\in H^{\infty}(X)$ satisfying (\ref{e1}) 
there are solutions $g_{1},\dots, g_{n}\in H^{\infty}(X)$ of
(\ref{e2}) with the bound $\max_{1\leq j\leq n}||g_{j}||_{\infty}\leq 
C(Y,n,\delta)$.
\end{Th}

This result, in a sense, completes our work started in [Br1]-[Br3] on the
corona problems on coverings of certain Riemann surfaces. Similarly
to [Br1]-[Br3] the methods used in 
the present paper are based on $L^{2}$ cohomology technique on complete 
K\"{a}hler manifolds and Cartan's $A$ and $B$ theorems for coherent 
Banach sheaves on Stein manifolds.
\begin{R}\label{re1}
{\rm
(1) Note that the assumption of the Caratheodory hyperbolicity of $Y$ cannot be
removed: It follows from the results of L\'{a}russon [L]
and the author [Br3] that for any integer $n\geq 2$ there are a compact
Riemann surface $S_{n}$ and its regular covering $r_{n}:\widetilde S_{n}\to
S_{n}$ such that
\begin{itemize}
\item[(a)]
$\widetilde S_{n}$ is a complex submanifold of an open Euclidean ball 
$\Bo_{n}\subset\Co^{n}$;
\item[(b)]
the embedding $i:\widetilde S_{n}\hookrightarrow\Bo_{n}$ induces an isometry
$i^{*}:H^{\infty}(\Bo_{n})\to H^{\infty}(\widetilde S_{n})$.
\end{itemize}
In particular, (b) implies that the maximal ideal spaces of 
$H^{\infty}(\widetilde S_{n})$ and $H^{\infty}(\Bo_{n})$ coincide.
Thus the corona problem is not solvable for $H^{\infty}(\widetilde S_{n})$.

(2) Under the assumptions of the theorem, let $U\hookrightarrow X$ be a domain 
such that the embedding induces an injective homomorphism of the 
corresponding fundamental groups and $r(U)\subset\subset Y$. Then as 
was shown in [Br1] and [Br2] in this case
the following extension of Theorem \ref{te1} is valid.
\begin{Th}\label{te2}
Let $A=(a_{ij})$ be an $n\times k$ matrix, $k<n$, with entries in 
$H^{\infty}(U)$. Assume that the family of minors of order $k$ of $A$
satisfies the corona condition. Then there is an $n\times n$
matrix $\widetilde A=(\widetilde a_{ij})$, 
$\widetilde a_{ij}\in H^{\infty}(U)$, so that $\widetilde a_{ij}=a_{ij}$ for
$1\leq j\leq k$, $1\leq i\leq n$, and $det\ \!\widetilde A=1$.
\end{Th} 

The proof of the theorem is based on a Forelli type theorem on projections
in $H^{\infty}$ discovered in [Br1] and a Grauert type theorem for
``holomorphic`` vector bundles on maximal ideal spaces (which are not usual
manifolds) of certain Banach algebras proved in [Br2]. In a forthcoming paper
we prove a result similar to Theorem \ref{te2} for matrices with entries
in $H^{\infty}(X)$ with $X$ satisfying assumptions of Theorem \ref{te1}.
The techniques are necessarily more complicated than those used in this paper.

(3) The remarkable class of Riemann surfaces $X$ for which a Forelli
type theorem is valid was introduced by Jones and Marshall [JM]. The definition
is in terms of an interpolating property for the critical points of the
Green function on $X$. For such $X$ the corona problem is solvable, as well.
Moreover, every $X$ from this class is of {\em Widom type}, see [W]
for the corresponding definition. (Roughly speaking, this means
that the topology of $X$ grows slowly as measured by the Green function.)
It is an interesting open question whether the surfaces $X$ in 
Theorem \ref{te1} are also of Widom type.

(4) Similarly to [JM] and [Br1] our proof of Theorem \ref{te1} uses the 
Carleson Corona Theorem for the open unit disk. 
}
\end{R}

{\bf 1.2.} In this part we formulate some results used in the proof
of Theorem \ref{te1}. First, we recall the following
\begin{D}\label{d1}
Let $X$ be a complex manifold. A sequence $\{x_{j}\}_{j\in\N}\subset X$ is
called interpolating for $H^{\infty}(X)$ if for every bounded sequence of
complex numbers $a=\{a_{j}\}_{j\in\N}$ there is an $f\in H^{\infty}(X)$ so that
$f(x_{j})=a_{j}$ for all $j$. The constant of interpolation for 
$\{x_{j}\}_{j\in\N}$ is defined as
\begin{equation}\label{e2'}
\sup_{||a||_{l^{\infty}}\leq 1}\inf\{||f||_{\infty}\ :\ f\in H^{\infty}(X),\ 
f(x_{j})=a_{j},\ j\in\N\}
\end{equation}
where $||a||_{l^{\infty}}:=\sup_{j\in\N}|a_{j}|$.
\end{D}

Let $r:X\to Y$ be an unbranched covering of a Caratheodory hyperbolic Riemann
surface of finite type $Y$. Let $K\subset\subset Y$ be a compact subset.
\begin{Th}\label{te3} 
For every $x\in K$ the sequence $r^{-1}(x)\subset X$ is interpolating for
$H^{\infty}(X)$ with the constant of interpolation bounded by a number
depending on $K$ and $Y$ only.
\end{Th}

To formulate our next result used in the proof
we will assume that $Y$ is equipped with a hermitian metric $h_{Y}$ with
the associated $(1,1)$-form $\omega_{Y}$. Then we equip $X$ with the
hermitian metric $h_{X}$ induced by the pullback 
$r^{*}\omega_{Y}$ of $\omega_{Y}$ to $X$. 
Now, if $\eta$ is a smooth differential $(0,1)$-form on $X$, by
$|\eta|_{z}$, $z\in X$, we denote the norm of $\eta$ at $z$
defined by the hermitian metric $h_{X}^{*}$ on the fibres of the cotangent
bundle $T^{*}X$ on $X$. We say that $\eta$ is bounded if
\begin{equation}\label{e3}
||\eta||:=\sup_{z\in X}|\eta|_{z}<\infty.
\end{equation}
\begin{Th}\label{te4}
Let $\eta$ be a smooth bounded $(0,1)$-form on $X$ with support $supp\ \!\eta$
satisfying $r(supp\ \!\eta)\subset K$
for some compact subset $K\subset\subset Y$. Then the equation
$\overline\partial f=\eta$
has a smooth bounded solution $f$ on $X$ such that
\begin{equation}\label{e4}
||f||_{L^{\infty}}:=\sup_{z\in X}|f(z)|\leq C||\eta||
\end{equation}
with $C$ depending on $K$, $Y$ and $h_{Y}$ only.
\end{Th}
\begin{R}\label{re2}
{\rm The methods used in the proofs of Theorems \ref{te3}, \ref{te4} can be
applied also to prove similar results for holomorphic $L^{p}$-functions on
unbranched coverings of certain Stein manifolds. We present these results
in a forthcoming paper.}
\end{R}
\sect{\hspace*{-1em}. Auxiliary Results.}
\quad In this part we collect some auxiliary results used in the proofs.

{\bf 2.1.} Let $Y$ be a Caratheodory hyperbolic Riemann surface of finite type.
According to the theorem of Stout [St, Th.$\!$ 8.1], there exist a compact
Riemann surface $R$ and a holomorphic embedding $\phi:Y\to R$ such that
$R\setminus\phi(Y)$ consists of finitely many closed disks with analytic
boundaries together with finitely many isolated points. Since $Y$ is
Caratheodory hyperbolic, the set of the disks in $R\setminus\phi(Y)$ is not
empty. Also, without loss of generality we may and will assume that
the set of isolated points in $R\setminus\phi(Y)$ is not empty, as well.
(For otherwise, $\phi(Y)$ is a bordered Riemann surface and the required
results follow from [Br1].) We will naturally identify
$Y$ with $\phi(Y)$. Also, we set 
\begin{equation}\label{eq2.1}
R\setminus Y:=\left(\bigsqcup_{1\leq i\leq k}\overline{D}_{i}\right)\bigcup
\left(\bigcup_{1\leq j\leq l}\{x_{j}\}\right)\ \ \ {\rm and}\ \ \
\widetilde Y:=Y\bigcup \left(\bigcup_{1\leq j\leq l}\{x_{j}\}\right)
\end{equation}
where each $D_{i}$ is biholomorphic to the open unit disk $\Di\in\Co$
and these biholomorphisms are extended to diffeomorphisms of the closures
$\overline{D}_{i}\to\overline\Di$. 

According to these definitions, $\widetilde Y$ is a bordered Riemann surface,
and there
is a bordered Riemann surface $\widetilde Y_{1}\subset\subset R$ with 
$\widetilde Y\subset\subset\widetilde Y_{1}$ and 
$\pi_{1}(\widetilde Y_{1})\cong\pi_{1}(\widetilde Y)$.

{\bf 2.2.} Next, we introduce a complete K\"{a}hler metric on 
$Y_{1}:=\widetilde Y_{1}\setminus (\cup_{1\leq j\leq l}\{x_{j}\})$. 

To this end we consider an open cover ${\cal U}=(U_{j})_{0\leq j\leq l}$ of 
$R$ such that for $1\leq j\leq l$ 
$U_{j}\subset\subset\widetilde Y$ is an open coordinate disk 
centered at $x_{j}$  and $U_{0}$ is a 
bordered Riemann surface which intersects each $U_{j}$, $1\leq j\leq l$, 
by a set biholomorphic to an open annulus. Let $\{\rho_{j}\}_{0\leq j\leq l}$
be a smooth partition of unity on $R$ subordinate to the cover ${\cal U}$.
By $z$ we denote a complex coordinate in $U_{j}$, $1\leq j\leq l$, such that
$z(x_{j})=0$ and $|z|<1$. 
We set $f_{j}:=\rho_{j}|z|^{2}$, $1\leq j\leq l$. Then
$f_{j}$ is a smooth nonnegative function on $R$ with 
$supp\ \!f_{j}\subset U_{j}$.

Now, we consider the positive smooth function 
\begin{equation}\label{e2.1}
f:=\frac{\rho_{0}}{2}+\sum_{j=1}^{l}f_{j}
\end{equation}
on $R$ and determine the $(1,1)$-form $\widetilde\omega$ on $R$ by
\begin{equation}\label{e2.2}
\widetilde\omega:=-\frac{\sqrt{-1}}{2\pi}
\partial\overline{\partial}\log(\log f)^{2}.
\end{equation}
Since by the definition $0<f<1$, the form $\omega$ is well defined. Also, in
an open neighbourhood of $x_{j}$, $1\leq j\leq l$, the form
$\omega$ is equal to
\begin{equation}\label{e2.3}
\omega_{P}:=
\frac{\sqrt{-1}}{\pi}\frac{dz\wedge d\overline{z}}{|z|^{2}(\log |z|^{2})^{2}}
=-\frac{\sqrt{-1}}{2\pi}\partial\overline\partial\log(\log |z|^{2})^{2}.
\end{equation}
In the natural identification 
$U_{j}\setminus\{x_{j}\}=\Di\setminus\{0\}$, $\omega_{P}$
coincides with the $(1,1)$-form of the 
Poincar\'{e}  metric on the punctured disk. 

Let $\omega_{R}$ be a K\"{a}hler $(1,1)$-form on the compact Riemann surface
$R$ from section 2.1. Since $\widetilde Y_{1}$ is a bordered Riemann surface
in $R$, there is a smooth
plurisubharmonic function $f_{R}$ defined in a neighbourhood
of the closure of $\widetilde Y_{1}$ such that
\begin{equation}\label{e2.4} 
\omega_{R}=\sqrt{-1}\cdot\partial\overline\partial f_{R}\ \ \ {\rm on}\ \ \ 
\widetilde Y_{1}.
\end{equation}

Further, since $\widetilde Y_{1}$ is a Stein manifold, by the Narasimhan 
theorem [N] there is a holomorphic embedding 
$i:\widetilde Y_{1}\hookrightarrow\Co^{3}$ of $\widetilde Y_{1}$ as 
a closed complex submanifold of $\Co^{3}$. By $\omega_{e}$ we denote the
$(1,1)$-form on $\widetilde Y_{1}$ obtained as the pullback by $i$ of the
Euclidean K\"{a}hler form $\sqrt{-1}(dz_{1}\wedge d\overline{z}_{1}+
dz_{2}\wedge d\overline{z}_{2}+dz_{3}\wedge d\overline{z}_{3})$ on $\Co^{3}$
(here $z_{1},z_{2},z_{3}$ are complex coordinates on $\Co^{3}$). Clearly,
$\omega_{e}$ is a K\"{a}hler form on $\widetilde Y_{1}$ (i.e., it is
positive on $\widetilde Y_{1}$ and $d$-closed).
\begin{Lm}\label{le1}
There is a positive number $c_{1}$ depending on $Y_{1}$, $\omega_{R}$ and
$\widetilde\omega$ such that the $(1,1)$-form
\begin{equation}\label{e2.5}
\omega:=\widetilde\omega+c_{1}(\omega_{R}+\omega_{e})
\end{equation}
is a complete K\"{a}hler form on $Y_{1}$.
\end{Lm}
{\bf Proof.} Since $\omega_{R}$ is a K\"{a}hler form on $R$,
by the definition of $\widetilde\omega$ there is a constant $c_{1}>0$ 
depending on $Y_{1}$, $\omega_{R}$ and
$\widetilde\omega$ such that
$$
\widetilde\omega>-c_{1}\omega_{R}\ \ \ {\rm on}\ \ \ Y_{1}.
$$
Thus the form $\omega=\widetilde\omega+c_{1}(\omega_{R}+\omega_{e})$ 
is positive (and $d$-closed) on $Y_{1}$. Its completeness
means that the path metric $d$ on $Y_{1}$ induced by $\omega$ is complete.
Let us check this fact.

Assume, on the contrary, that $d$ is not complete. This means that
there is a sequence $\{w_{n}\}\subset Y_{1}$
convergent either to the boundary of $\widetilde Y_{1}$ or to one of the
points $x_{j}$, $1\leq j\leq l$, such that $\{d(o,w_{n})\}$ is bounded
(for a fixed point $o\in Y_{1}$). Then, since $\omega\geq\omega_{e}$, the
sequence $\{i(w_{n})\}\subset\Co^{3}$ is bounded. This implies that
$\{w_{n}\}$ cannot converge to the boundary of $\widetilde Y_{1}$. Thus it
converges to one of $x_{j}$. But since $\omega\geq\omega_{P}$ near $x_{j}$,
the latter is impossible because the Poincar\'{e} metric on the punctured
disk is complete.\ \ \ \ \ $\Box$

{\bf 2.3.} According to our construction the embedding $Y\hookrightarrow Y_{1}$
induces an isomorphism of the corresponding fundamental groups. By the 
covering homotopy theorem
this implies that for any unbranched covering $r:X\to Y$ there is an
unbranched covering $r_{1}:X_{1}\to Y_{1}$ and an embedding 
$j:X\hookrightarrow X_{1}$ such that $r_{1}\circ j=r$ and 
$j_{*}:\pi_{1}(X)\to\pi_{1}(X_{1})$ is an isomorphism. Without loss of
generality we identify $j(X)$ with $X$. Then $r:=r_{1}|_{X}$. Now the form
$r_{1}^{*}\omega$ with $\omega$ determined by (\ref{e2.5}) is a complete
K\"{a}hler form on $X_{1}$.

Let $TY_{1}$ be the complex tangent bundle on $Y_{1}$ equipped with
the hermitian metric induced by the K\"{a}hler form $\omega_{R}$. 
Since $TY_{1}$ is the restriction to $Y_{1}$ of the tangent bundle
$TR$ on $R$ with the hermitian metric defined by $\omega_{R}$, the
curvature $\Theta_{Y_{1}}$ of $TY_{1}$ satisfies
\begin{equation}\label{e2.6}
\Theta_{Y_{1}}\geq -c_{2}\omega_{R}
\end{equation}
for some $c_{2}>0$ depending on $Y_{1}$ and $\omega_{R}$. In turn, the
curvature $\Theta_{X_{1}}:=r_{1}^{*}\Theta_{Y_{1}}$ of the tangent bundle
$TX_{1}$ on $X_{1}$ equipped with the hermitian metric induced by 
$r_{1}^{*}\omega_{R}$ satisfies
\begin{equation}\label{e2.7}
\Theta_{X_{1}}\geq -c_{2}r_{1}^{*}\omega_{R}.
\end{equation}

Next, by the definition of $\omega_{e}$ there is a smooth
plurisubharmonic function $g$ on $\widetilde Y_{1}$ such that
$\omega_{e}=\sqrt{-1}\cdot\partial\overline{\partial} g$ (as such $g$ one 
takes the pullback by $i$ of the function $||z||^{2}:=|z_{1}|^{2}+
|z_{2}|^{2}+|z_{3}|^{2}$ on $\Co^{3}$). 

Let $E_{0}:=X_{1}\times\Co$ be the
trivial holomorphic line bundle on $X_{1}$. We equip $E_{0}$ with the hermitian
metric $e^{h_{1}-(c_{1}+c_{2})g_{1}}$ where 
$h_{1}=r_{1}^{*}h:=\frac{\log(\log r_{1}^{*}f)^{2}}{2\pi}$
with $f$ from (\ref{e2.1}) and $g_{1}:=r_{1}^{*}(g+f_{R})$, see (\ref{e2.4}). 
(This means that
for $z\times v\in E$ the square of its norm in this metric equals
$e^{h_{1}(z)-(c_{1}+c_{2})g_{1}(z)}|v|^{2}$ 
where $|v|$ is the modulus of $v\in\Co$.)
Then by (\ref{e2.7})
the curvature $\Theta_{E}$ of the bundle $E:=E_{0}\otimes TX_{1}$ 
satisfies
\begin{equation}\label{e2.8}
\Theta_{E}:=-\sqrt{-1}\cdot\partial\overline\partial
\log e^{h_{1}-(c_{1}+c_{2})g_{1}}+\Theta_{X_{1}}\geq
r_{1}^{*}\omega.
\end{equation}

{\bf 2.4.} Let $X$ be a complete K\"{a}hler manifold of dimension $n$ with a 
K\"{a}hler form $\omega$ and
$E$ be a hermitian holomorphic vector bundle on $X$ with curvature $\Theta$. 
Let $L_{2}^{p,q}(X,E)$ be the space of $L^{2}$ $E$-valued $(p,q)$-forms on $X$
with the $L^{2}$ norm, and let $W_{2}^{p,q}(X,E)$ be the subspace of forms
such that $\overline\partial\eta$ is $L^{2}$. (The forms $\eta$ may be taken 
to be either smooth or just measurable, in which case 
$\overline\partial\eta$ is understood in the distributional sense.) The
cohomology of the resulting $L^{2}$ Dolbeault complex 
$(W_{2}^{\cdot ,\cdot},\overline\partial)$ is the $L^{2}$ cohomology
$$
H_{(2)}^{p,q}(X,E)=Z_{2}^{p,q}(X,E)/B_{2}^{p,q}(X,E),
$$
where $Z_{2}^{p,q}(X,E)$ and $B_{2}^{p,q}(X,E)$ are the spaces of 
$\overline\partial$-closed and $\overline\partial$-exact forms in
$L_{2}^{p,q}(X,E)$, respectively.

If $\Theta\geq\epsilon\omega$
for some $\epsilon>0$ in the sense of Nakano, then the $L^{2}$ 
Kodaira-Nakano vanishing theorem, see [D], [O], states that
\begin{equation}\label{e2.9}
H_{(2)}^{n,r}(X,E)=0\ \ \ {\rm for}\ \ \ r>0 .
\end{equation}

Moreover, for $\eta\in Z_{2}^{n,r}(X,E)$, $r>0$, there is a form
$\widetilde\eta\in W_{2}^{n,r-1}(X,E)$ such that 
$\overline\partial\widetilde\eta=\eta$ and
\begin{equation}\label{e2.10}
||\widetilde\eta||_{2}\leq\frac{1}{\epsilon}||\eta||_{2},
\end{equation}
see [D, Rem.$\!$ 4.2]. Here symbols $||\cdot||_{2}$ denote
the corresponding $L^{2}$ norms.

We can apply this result to the bundle $E$ from section 2.3 with $X_{1}$
equipped with the complete K\"{a}hler form $r_{1}^{*}\omega$. Then 
from (\ref{e2.8}) we obtain
\begin{Proposition}\label{pr2.2}
For every $\eta\in W_{2}^{1,1}(X_{1},E)$ there is $\widetilde\eta\in 
W_{2}^{1,0}(X_{1},E)$ such that $\overline\partial\widetilde\eta=\eta$ and
\begin{equation}\label{e2.11}
||\widetilde\eta||_{2}\leq ||\eta||_{2}.
\end{equation}
\end{Proposition}

{\bf 2.5.}
Let $T^{*}X_{1}$ be the cotangent bundle on $X_{1}$ equipped with the 
hermitian metric induced by $r_{1}^{*}\omega$. We consider the hermitian line
bundle $V:=E\otimes T^{*}X_{1}$ equipped with the tensor product of the
corresponding hermitian metrics. Then from Proposition \ref{pr2.2} we obtain
\begin{equation}\label{e2.12}
H_{(2)}^{0,1}(X_{1},V)\cong H_{(2)}^{1,1}(X_{1},E)=0
\end{equation}

Moreover, for every $\eta\in W_{2}^{0,1}(X_{1},V)$ there is 
$F\in W_{2}^{0,0}(X_{1},V)$ such that $\overline\partial F=\eta$ and
\begin{equation}\label{e2.13}
||F||_{2}\leq ||\eta||_{2}.
\end{equation}

Further, there is a canonical isomorphism $I:X_{1}\times\Co\to V$
defined in local coordinates $z$ on $X_{1}$ by the formula
\begin{equation}\label{e2.14}
I(z\times v):=v\cdot 1\otimes\frac{\partial}{\partial z}\otimes dz.
\end{equation}
(Clearly this definition does not depend on the choice of local coordinates.)
In what follows we identify $V$ with $X_{1}\times\Co$ by $I$.
\sect{\hspace*{-1em}. Proof of Theorem \ref{te3}.}
\quad Below by $A, B, C, c$ etc we denote constants depending 
on characteristics related to the sets $Y$, $x\in Y$ and 
$K\subset\subset Y$ but not on coverings
$X$ of $Y$. (We will briefly say that they depend on $Y$, $x$ and 
$K$ only.) These constants may change from line to line and even in a
single line.

{\bf 3.1.} Let $x\in K\subset\subset Y$. 
We must check that the sequence
$r^{-1}(x)\subset X$ is interpolating for $H^{\infty}(X)$ with the constant
of interpolation bounded by a number depending on $K$ and $Y$ only.
Fix a neighbourhood $\widehat Y$ of the closure of $Y$ in $Y_{1}$ such that
$\widehat Y$ is relatively compact in $\widetilde Y_{1}$, see section 2.1 for
the corresponding definitions.
First we will prove that $r^{-1}(x)$ is interpolating for the space
$L_{{\cal O}}^{2}(\widehat X_{1};r_{1}^{*}\omega_{R})$ of holomorphic 
$L^{2}$-functions on the covering $\widehat X_{1}:=r_{1}^{-1}(\widehat Y)$ of 
$\widehat Y$ with norm defined by integration with respect to 
$r_{1}^{*}\omega_{R}$ (recall that $X\subset X_{1}$ and
$r_{1}|_{X}=r$).
\begin{Proposition}\label{pr3.1}
Let $a$ be an $l^{2}$-function on $r^{-1}(x)$ with norm $||a||_{l^{2}}$.
Then there is a function 
$f\in L_{{\cal O}}^{2}(\widehat X_{1};r_{1}^{*}\omega_{R})$ 
such that $f|_{r^{-1}(x)}=a$ and
$$
||f||_{2}\leq c||a||_{l^{2}}
$$
with $c$ depending on $K$ and $Y$ only.
\end{Proposition}
{\bf Proof.} Using the fact that 
the closure of $\widetilde Y_{1}$ in $R$ possesses
a Stein neighbourhood and applying some basic results of
the theory of Stein manifolds (see, e.g., [GR]) one obtains easily that
there are a holomorphic function $\phi_{x}$ with
a simple zero at $x$ defined in a fixed
neighbourhood of the closure of $\widetilde Y_{1}$,
a simply connected coordinate neighbourhood $U_{x}\subset\subset Y$
of $x$ with a complex coordinate $z$, $z(x)=0$ and
$|z|<1$ on $U_{x}$, and positive numbers $A$ and $r$,
$0<r<1$, depending on $x$ and $Y$ only such that
\begin{itemize}
\item[(1)]
$$
\sup_{y\in\widetilde Y_{1}}|\phi_{x}(y)|\leq A;
$$
\item[(2)]
$$
|\phi_{x}(y)|\geq\frac{1}{A}\ \ \ {\rm for\ all}\ \ \ y\in U_{x;r}:=
\{z\in U_{x}\ :\ |z|\geq r\};
$$
\item[(3)]
$$
\phi_{x}(y)\neq 0\ \ \ {\rm for\ all}\ \ \ y\in U_{x}\setminus\{x\}\ \ \
$$ 

\end{itemize}
(For a construction of such $\phi_{x}$ see, e.g., [Br2, Cor.$\!$ 1.8].)

Further, there is a $C^{\infty}$-function $\rho_{x}$,
$0\leq\rho_{x}\leq 1$, on $U_{x}$ such that 
$\rho_{x}$ is equal to 1 on $U_{x}\setminus U_{x;r}$ and 0 outside 
$U_{x}$, and
\begin{equation}\label{e3.1}
|d\rho_{x}|_{z;\omega}\leq B\ \ \ {\rm for\ all}\ \ \ z\in Y_{1}
\end{equation}
where $\{|\cdot|_{z;\omega}\ :\ z\in Y_{1}\}$ is the hermitian metric on the fibres 
of the
cotangent bundle $T^{*}Y_{1}$ on $Y_{1}$
determined by the form $\omega$ (see (\ref{e2.5})), and the constant
$B$ depends on $x$ and $Y$ only.

Since $U_{x}$ is simply connected, $r_{1}^{-1}(U_{x})$ is biholomorphic to
$U_{x}\times S$ where $S$ is the fibre of $r_{1}$. In what follows without 
loss of generality we will identify these sets. Then $a$ is an
$l^{2}$-function of $\{x\}\times S$. 
We extend $a$ to a locally constant function 
$\hat a$ on $U_{x}\times S$ by the formula
\begin{equation}\label{e3.2}
\hat a(z,s):=a(x,s)\ \ \ {\rm for\ all}\ \ \  (z,s)\in U_{x}\times S.
\end{equation}

Let us consider a $(0,1)$ form $\hat\eta$ on $X_{1}$ determined by 
\begin{equation}\label{e3.3}
\hat\eta(w)=\left\{
\begin{array}{ccc}
\displaystyle
\frac{\hat a(z,s) d\rho_{x}(z)}{\phi_{x}(z)}&{\rm if}&
w=(z,s)\in U_{x}\times S\\
\\
\displaystyle
0&{\rm if}&w\not\in U_{x}\times S.
\end{array}
\right.
\end{equation}

Next, we fix a noncompact neighbourhood $O$ of the closure $\widetilde Y_{1}$
in $R$. Since $O$ is a one-dimensional Stein manifold, it is homotopically
equivalent to a one-dimensional CW-complex (see, e.g., [GR]). In particular,
any continuous vector bundle on $O$ is topologically trivial. Then by the
Grauert theorem [Gr], any holomorphic vector bundle on $O$ is also trivial.
Applying this to the bundle $T^{*}R|_{O}$ we find a nowhere vanishing 
holomorphic section $\lambda$ of $T^{*}R|_{O}$. (By the definition 
$\lambda$ is a holomorphic $1$-form on $O$.) Then $\lambda^{-1}=
\frac{1}{\lambda}$ is
a nowhere vanishing holomorphic section of the tangent bundle $TO$ on $O$.
Moreover, there is a positive constant $C$ depending on $\widetilde Y_{1}$ 
such that
\begin{equation}\label{e3.4}
\frac{1}{C}\leq |\lambda|_{z;\omega_{R}}\leq C,\ \ \ \ \  
\frac{1}{C}\leq |\lambda^{-1}|_{z;\omega_{R}}\leq C\ \ \ 
{\rm for\ all}\ \ \ z\in
\widetilde Y_{1}
\end{equation}
where symbols 
$|\cdot|_{z;\omega_{R}}$ denote the corresponding norms on 
$T^{*}\widetilde Y_{1}$
and $T\widetilde Y_{1}$ at $z\in\widetilde Y_{1}$ determined by the form
$\omega_{R}$. 

Continuing the proof of the proposition consider the $V$-valued
$(0,1)$-form $\eta:=\hat\eta\wedge r_{1}^{*}\lambda\otimes 
r_{1}^{*}\lambda^{-1}$ on $X_{1}$.
By $N_{w}$ and $N_{w}'$ we denote the hermitian norms on $V$ and $E$
(see section 2.3) at $w\in X_{1}$.
Then we have
\begin{equation}\label{e3.5}
\begin{array}{cc}
\displaystyle
N_{w}(\eta):=N_{w}'(\hat\eta\otimes r_{1}^{*}\lambda^{-1})\cdot 
|r_{1}^{*}\lambda|_{w;r_{1}^{*}\omega}=\\
\\
\displaystyle
|\hat\eta|_{w;r_{1}^{*}\omega}\cdot 
|r_{1}^{*}\lambda^{-1}|_{w;r_{1}^{*}\omega_{R}}
\cdot |r_{1}^{*}\lambda|_{w;r_{1}^{*}\omega}
\cdot e^{h_{1}(w)-(c_{1}+c_{2})g_{1}(w)}.
\end{array}
\end{equation}
Here $|\hat\eta|_{w;r_{1}^{*}\omega}$ and 
$|r_{1}^{*}\lambda|_{w;r_{1}^{*}\omega}$ are determined
by the form $r_{1}^{*}\omega$ and 
$|r_{1}^{*}\lambda^{-1}|_{w;r_{1}^{*}\omega_{R}}$ is determined by
the form $r_{1}^{*}\omega_{R}$. Since $\omega_{R}$ is equivalent to
$\omega$ on $U_{x}\subset\subset Y$, from (\ref{e3.4}) we
obtain  for all $w\in r_{1}^{-1}(U_{x})$:
\begin{equation}\label{e3.6}
|r_{1}^{*}\lambda|_{w;r_{1}^{*}\omega}:=|\lambda|_{r_{1}(w);\omega}\leq C_{1}
\ \ \ {\rm and}\ \ \ 
|r_{1}^{*}\lambda^{-1}|_{w;r_{1}^{*}\omega_{R}}:=
|\lambda^{-1}|_{r_{1}(w);\omega_{R}}\leq C
\end{equation}
where $C_{1}$ depends on $C$ and on the constant of the equivalence of
$\omega_{R}$ and $\omega$ on $U_{x}$.
Also, by the definition of $h_{1}$ and $g_{1}$, see section 2.3,
there is a constant $C_{2}>0$ depending on $h_{1}$, $g_{1}$,
$U_{x}$ and $Y_{1}$ such that
\begin{equation}\label{e3.7}
 e^{h_{1}(w)-(c_{1}+c_{2})g_{1}(w)}\leq C_{2}\ \ \ {\rm for\ all}\ \ \
w\in r_{1}^{-1}(U_{x}).
\end{equation}
From here, (\ref{e3.6}), (\ref{e3.5}),
(\ref{e3.1})-(\ref{e3.3}) and the definition of $\phi_{x}$ we obtain
\begin{equation}\label{e3.8}
N_{w}(\eta)\leq\widetilde C|\tilde a(w)|.
\end{equation}
where $\widetilde C:=ABCC_{1}C_{2}$,
$\tilde a:=\hat a\cdot\chi_{r_{1}^{-1}(U_{x})}$ and 
$\chi_{r_{1}^{-1}(U_{x})}$
is the characteristic function of $r_{1}^{-1}(U_{x})$.
Using (\ref{e3.8}) we estimate the $L^{2}$ norm of $\eta$ as follows:
\begin{equation}\label{e3.9}
||\eta||_{2}:=\left(\int_{w\in X_{1}}N_{w}^{2}(\eta)\ \!(r_{1}^{*}\omega)(w)
\right)
^{1/2}\leq\widetilde C||a||_{l^{2}}\left(\int_{y\in U_{x}}
\omega(y)\right)^{1/2}\leq \widehat C||a||_{l^{2}}.
\end{equation}
Here $\widehat C$ depends on $x$ and $Y$ only. 

Next, according to (\ref{e2.13}) there
is $F\in W_{2}^{0,0}(X_{1},V)$ such that $\overline\partial F=\eta$
and
\begin{equation}\label{e3.10}
||F||_{2}\leq |\eta||_{2}\leq \widehat C||a||_{l^{2}}.
\end{equation}
We regard $F$ as a function on $X_{1}$, see section 2.5. 

Observe that since $\widehat Y (\ \subset Y_{1})$ is relatively compact
in $\widetilde Y_{1}$,
$g_{1}$ is bounded on $\widehat X_{1}$ and $h_{1}$ is bounded from
below on $\widehat X_{1}$, see section 2 for the corresponding 
definitions. Moreover, $\omega_{R}\leq 
\widetilde c\cdot\lambda\wedge\overline\lambda$ on $\widetilde Y_{1}$
for some $\widetilde c$ depending on $\widetilde Y_{1}$.
These facts, (\ref{e3.4}) and (\ref{e3.5}) imply
\begin{equation}\label{e3.11}
\begin{array}{c}
\displaystyle
\left(\int_{\widehat X_{1}}|F|^{2}\! \ \cdot
r_{1}^{*}\omega_{R}\right)^{1/2}\leq C\left(\int_{w\in\widehat X_{1}}
[N_{w}'(F\otimes r_{1}^{*}\lambda^{-1})]^{2}\!\
(r_{1}^{*}\lambda\wedge r_{1}^{*}\overline\lambda)(w)\right)^{1/2}=
\\
\\
\displaystyle
C\left(\int_{w\in\widehat X_{1}}N_{w}^{2}(F\cdot r_{1}^{*}\lambda\otimes
r_{1}^{*}\lambda^{-1})\! \ 
(r_{1}^{*}\omega)(w)\right)^{1/2}\leq\\
\\
\displaystyle
C\left(\int_{w\in X_{1}}N_{w}^{2}(F\cdot r_{1}^{*}\lambda\otimes
r_{1}^{*}\lambda^{-1})
\! \ (r_{1}^{*}\omega)(w)\right)^{1/2}:=
C||F||_{2}
\end{array}
\end{equation}
with $C$ depending on $Y$. Let us consider the function
$f:=\widetilde a\cdot\rho_{x}-F\cdot r_{1}^{*}\phi_{x}$
on $\widehat X_{1}$. Then according to (\ref{e3.10}) and (\ref{e3.11})
$$
\left(\int_{\widehat X_{1}}|f|^{2}\! \ r_{1}^{*}\omega_{R}\right)^{1/2}\leq 
c||a||_{l^{2}}
$$
where $c$ depends on $Y$ and $x$ only. Moreover, from
(\ref{e3.3}) it follows that $f$ is holomorphic and $f|_{r_{1}^{-1}(x)}=a$.

This shows that $r^{-1}(x)$ is interpolating for 
$L_{{\cal O}}^{2}(\widehat X_{1}; r_{1}^{*}\omega_{R})$.

To complete the proof of the proposition it remains to prove that for
$x\in K\subset\subset Y$ the 
constant of interpolation of $r^{-1}(x)$ with respect to 
$L_{{\cal O}}^{2}(\widehat X_{1}; r_{1}^{*}\omega_{R})$ is bounded by a number 
depending on $K$ and $Y$ only.

Let us consider the restriction map $R_{x}: 
L_{{\cal O}}^{2}(\widehat X_{1}; r_{1}^{*}\omega_{R})\to
L_{{\cal O}}(\widehat X_{1}; r_{1}^{*}\omega_{R})|_{r^{-1}(x)}$.
\begin{Lm}\label{le3.2}
$R_{x}$ maps $L_{{\cal O}}^{2}(\widehat X_{1}; r_{1}^{*}\omega_{R})$ 
continuously onto $l^{2}(r^{-1}(x))$.
Moreover, the norm $||R_{x}||$ of $R_{x}$ is bounded by a constant depending
on $Y$, $x$.
\end{Lm}
{\bf Proof.}
We will consider the coordinate neighbourhood $U_{x}\subset\subset Y$ from
the proof of Proposition \ref{pr3.1} with a complex coordinate $z$ so that
$z(x)=0$, $|z|<1$ on $U_{x}$. Also, we naturally identify $r_{1}^{-1}(U_{x})$
with $U_{x}\times S$ where $S$ is the fibre of $r_{1}$. By definition there
is a constant $C$ depending on  $Y$, $U_{x}$ and $\omega_{R}$ such that
\begin{equation}\label{e3.12}
\sqrt{-1}\cdot dz\wedge d\overline{z}\leq C\omega_{R}\ \ \ {\rm on}\ \ \
U_{x}.
\end{equation}
Let $f\in L_{{\cal O}}^{2}(\widehat X_{1}; r_{1}^{*}\omega_{R})$. 
Then on $r^{-1}(x)=\{x\}\times S$ we have 
by the mean-value property for subharmonic functions
\begin{equation}\label{e3.13}
\begin{array}{c}
\displaystyle
\sum_{s\in S}|f(x,s)|^{2}\leq\sum_{s\in S}\left(\frac{1}{\pi}\int_{U_{x}}
|f(z,s)|^{2}\! \ \sqrt{-1}\cdot dz\wedge d\overline z\right)\leq\\
\\
\displaystyle
\frac{C}{\pi}\int_{U_{x}}
\left(\sum_{s\in S}|f(z,s)|^{2}\right)\omega_{R}(z)=\frac{C}{\pi}
\int_{r_{1}^{-1}(U_{x})}|f|^{2}\! \ r_{1}^{*}\omega_{R}\leq\frac{C}{\pi}
||f||_{2}^{2}.
\end{array}
\end{equation}

This shows that $R_{x}$ maps 
$L_{{\cal O}}^{2}(\widehat X_{1}; r_{1}^{*}\omega_{R})$ continuously into
$l^{2}(r^{-1}(x))$. Also, $R_{x}$ is surjective according to
the first part of Proposition \ref{pr3.1} proved above.\ \ \ \ \ $\Box$

Now, since $R_{x}: L_{{\cal O}}^{2}(\widehat X_{1}; r_{1}^{*}\omega_{R})
\to l^{2}(r^{-1}(x))$ is 
a linear continuous surjective map of Hilbert spaces, there is a linear
continuous map $T_{x}:l^{2}(r^{-1}(x))\to  
L_{{\cal O}}^{2}(\widehat X_{1}; r_{1}^{*}\omega_{R})$ such that
$R_{x}\circ T_{x}=id$. Let $\{e_{s}\}_{s\in S}$, $e_{s}(x,t)=0$ for
$t\neq s$ and $e_{s}(x,s)=1$, be an orthonormal basis of $l^{2}(r^{-1}(x))$.
We set
$$
h_{s}:=T_{x}(e_{s})\in L_{{\cal O}}^{2}(\widehat X_{1}; r_{1}^{*}\omega_{R}).
$$
Then for a sequence $a=\{a_{s}\}_{s\in S}\in l^{2}(S)$ we have
\begin{equation}\label{e3.14}
h_{a}:=\sum_{s\in S}a_{s}h_{s}\in 
L_{{\cal O}}^{2}(\widehat X_{1}; r_{1}^{*}\omega_{R})\ \ \ {\rm and}\ \ \
||h_{a}||_{2}\leq c||a||_{l^{2}}.
\end{equation}

Further, for each $y\in U_{x}$ by 
$L_{y}:l^{2}(r^{-1}(y))\to l^{2}(r^{-1}(x))$ we
denote the natural isomorphism that sends 
$a(y,s)\in l^{2}(r^{-1}(y))$ to $a(x,s)\in l^{2}(r^{-1}(x))$.
Let us consider the map $S_{y}:=R_{y}\circ T_{x}\circ L_{y}:l^{2}(r^{-1}(y))\to
l^{2}(r^{-1}(y))$ determined by the formula
\begin{equation}\label{e3.15}
[S_{y}(a)](y,t):=\sum_{s\in S}a_{s}h_{s}(y,t),
\ \ \ (y,t)\in r^{-1}(y).
\end{equation}
Here $a(y,t)=\sum_{s\in S}a_{s}e_{s}(y,t)$, $t\in S$,
and $e_{s}(y,\cdot)$ are determined similarly
to $e_{s}(x,\cdot)$. Identifying $a\in l^{2}(r^{-1}(y))$ with 
$\{a_{s}\}_{s\in S}\in l^{2}(S)$ we can regard, according to
(\ref{e3.15}),  $\{S_{y}\}_{y\in U_{x}}$ as a family of bounded linear 
operators $l^{2}(S)\to l^{2}(S)$ depending holomorphically on $y\in U_{x}$.
According to (\ref{e3.14}) and Lemma \ref{le3.2} there is a constant
$c'$ depending on $U_{x}$ and $Y$ such that
$$
||S_{y}||\leq c'\ \ \ {\rm for\ all}\ \ \ y\in U_{x}.
$$
Moreover, by our construction $S_{x}=I$ where $I: l^{2}(S)\to l^{2}(S)$ is the 
identity operator. Identifying $U_{x}$ with
$\Di$ by the coordinate $z$, we obtain by the Cauchy integral formula for
bounded holomorphic on $\Di$ functions:
$$
S_{z}:=I+\sum_{k=1}^{\infty}S_{k}z^{k}\ \ \ {\rm for\ some}\ \ \
S_{k}:l^{2}(S)\to l^{2}(S),\ \ \ ||S_{k}||\leq c'.
$$
In particular, for $|z|<\frac{1}{2c'+4}$ we have
$$
\left|\left|\sum_{k=1}^{\infty}S_{k}z^{k}\right|\right|\leq c'\frac{|z|}{1-|z|}<\frac{2}{3}.
$$
Thus for every $y\in U_{x}$, $|z(y)|<\frac{1}{2c'+4}$, the inverse operator 
$S_{y}^{-1}$ exists and its norm is bounded by $\frac{1}{1-2/3}=3$.

Finally we set 
\begin{equation}\label{comp}
\widehat T_{y}:=T_{x}\circ L_{y}\circ S_{y}^{-1},\ \ \ \
y\in \widehat U_{x}:=\left\{y\in U_{x}\ :\ |z(y)|<\frac{1}{2c'+4}\right\}.
\end{equation}
Then by the definition we have 
$$
R_{y}\circ\widehat T_{y}=id\ \ \ {\rm for\ all}\ \ \ y\in\widehat U_{x}.
$$
This shows that $\{\widehat T_{y}:l^{2}(r^{-1}(y))\to 
L_{{\cal O}}^{2}(\widehat X_{1}; r_{1}^{*}\omega_{R})\ :\ 
y\in \widehat U_{x}\}$ is a family of interpolation operators depending 
holomorphically on 
$y$ such that $||\widehat T_{y}||\leq 3c$.
Taking a finite open cover of $K\subset\subset Y$ by the sets 
$\widehat U_{x}$, $x\in K$, and considering on these sets
the interpolation operators $\widehat T_{y}$, $y\in\widehat U_{x}$, 
we obtain that for every $x\in K$ the constant
of interpolation of $r^{-1}(x)$ with respect to 
$L_{{\cal O}}^{2}(\widehat X_{1}; r_{1}^{*}\omega_{R})$
is bounded by a number depending on $K$ and $Y$ only.

This completes the proof of the proposition.\ \ \ \ \ $\Box$

{\bf 3.2.} Let us prove now that $r^{-1}(x)$ is interpolating for
$H^{\infty}(X)$ with the constant of interpolation bounded by a number
depending on $K$ and $Y$ only.

We will use the interpolation operators $\widehat T_{y}$, 
$y\in \widehat U_{x}\subset\subset U_{x}$, of the previous section.
As before we set 
$$
h_{s,y}:=\widehat T_{y}(e_{s}(y,\cdot))\in 
L_{{\cal O}}^{2}(\widehat X_{1}; r_{1}^{*}\omega_{R}).
$$
Then the family of functions
$\{h_{s,y}\ :\ s\in S, y\in\widehat U_{x}\}$ depends 
holomorphically on $y$.

Now for a sequence $a=\{a_{s}\}_{s\in S}\in l^{2}(S)$ we have
\begin{equation}\label{hay}
h_{a,y}:=\sum_{s\in S}a_{s}h_{s,y}\in 
L_{{\cal O}}^{2}(\widehat X_{1}; r_{1}^{*}\omega_{R})\ \ \ {\rm and}
\ \ \ ||h_{a,y}||_{2}\leq 3c||a||_{l^{2}}.
\end{equation}

From here and Lemma \ref{le3.2} it follows that
\begin{equation}\label{e3.16}
\sum_{z\in r^{-1}(w)}|h_{a,y}(z)|^{2}\leq 
c||a||_{l^{2}}^{2},\ \ \
\ \
w\in Y,\ y\in \widehat U_{x},
\end{equation}
with $c$ depending on $w$, $x$ and $Y$.

To continue the proof we require an extension of Lemma \ref{le3.2}.

For its formulation we fix a holomorphic function $\phi$ defined in
a neighbourhood of the closure of $\widetilde Y_{1}\subset\subset R$ 
having simple zeros at all points $x_{j}$, $1\leq j\leq l$, see (\ref{eq2.1}),
and nonzero outside these points. (Such $\phi$ exists, e.g., by
[Br2, Cor.$\!$ 1.8].) As before, by $R_{w}: 
L_{{\cal O}}^{2}(\widehat X_{1};r_{1}^{*}\omega_{R})\to l^{2}(r^{-1}(w))$,
$w\in Y$, we denote the restriction map (here $r:=r_{1}|_{X}$).
\begin{Lm}\label{le3.3}
There is a constant $A>0$ depending on $Y$ such that
$$
||R_{w}||\leq\frac{A}{|\phi(w)|}.
$$
\end{Lm}
{\bf Proof.} Let $U_{j}\subset\subset\widetilde Y$ be a coordinate 
neighbourhood of $x_{j}$
with a complex coordinate $z$ so that $z(x_{j})=0$, $|z|<1$ on $U_{j}$.
We have
$$
\sqrt{-1}\cdot dz\wedge d\overline{z}\leq C\omega_{R}\ \ \ {\rm on}\ \ \
U_{j}
$$
for some $C$ depending on $\omega_{R}$, $z$ and $Y$. Thus for any
$f\in L_{{\cal O}}^{2}(\widehat X_{1}; r_{1}^{*}\omega_{R})$ 
its restriction $f|_{r_{1}^{-1}(U_{j})}$ belongs
to the $L^{2}$ space on $r_{1}^{-1}(U_{j})$ 
defined by integration with respect to the form
$r_{1}^{*}(\sqrt{-1}\cdot dz\wedge d\overline{z})$, and the 
$L^{2}$ norm of the restriction is bounded by $C||f||_{2}$.

Next, for a point $w\in \widetilde U_{j}\setminus\{x_{j}\}$,
$\widetilde U_{j}:=\{z\in U_{j}\ :\ |z|<1/2\}$,
 we set $d:=|z(w)|$ and
$D_{w}=\{y\in U_{j}\ :\ |z(y)-z(w)|<d\}$.
Then $D_{w}\subset U_{j}\setminus\{x_{j}\}$ and so
$r_{1}^{-1}(D_{w})$ is naturally identified with $D_{w}\times S$. In this 
identification we have by the mean-value property for subharmonic
functions:
$$
\begin{array}{c}
\displaystyle
\sum_{s\in S}|f(w,s)|^{2}\leq\frac{1}{\pi d^{2}}\int_{D_{w}}
\left(\sum_{s\in S}|f(z,s)|^{2}\right)
\! \ \sqrt{-1}\cdot dz\wedge d\overline{z}\leq\\
\\
\displaystyle
\frac{C}{\pi d^{2}}\int_{r_{1}^{-1}(D_{w})}|f|^{2}\! \ r_{1}^{*}\omega_{R}
\leq\frac{C}{\pi d^{2}}||f||_{2}^{2}.
\end{array}
$$
From here and the fact that the function $\phi/z$ is bounded on $U_{j}$
we obtain that there is a constant $c_{j}>0$ depending on $Y$ such that
\begin{equation}\label{e3.17}
||f|_{r_{1}^{-1}(w)}||_{l^{2}}\leq\frac{c_{j}}{|\phi(w)|}||f||_{2}\ \ \ 
{\rm for\ all}\ \ \ w\in \widetilde U_{j}.
\end{equation}

This proves the required inequality for $w\in\cup_{1\leq j\leq l}\ \!
\widetilde U_{j}$.

The remaining part $K:=Y\setminus (\cup_{1\leq j\leq l}\ \!
\widetilde U_{j})$ is a relatively compact subset of 
$\widehat Y\setminus (\cup_{1\leq j\leq l}\{x_{j}\})$, see section 3.1
for the definition of $\widehat Y$.
Then the required estimate on $K$ follows from (\ref{e3.13}) and the 
fact that $|\phi|$ is bounded on $K$.
We leave the details to the reader.\ \ \ \ \ $\Box$

From this lemma we obtain the following improvement of (\ref{e3.16}):
\begin{equation}\label{e3.16'}
\sum_{z\in r^{-1}(w)}|h_{a,y}(z)|^{2}\leq 
\frac{C}{|\phi(w)|^{2}}||a||_{l^{2}}^{2},\ \ \
\ \
w\in Y,\ y\in \widehat U_{x},
\end{equation}
with $C$ depending on $Y$ and $x$ only. This and the definition of
$h_{a,y}$, see (\ref{hay}), imply
\begin{equation}\label{e3.16''}
\sum_{s\in S}|h_{s,y}(z)|^{2}\leq 
\frac{C}{|\phi(w)|^{2}},\ \ \
\ \
z\in r^{-1}(w),\ w\in Y,\ y\in \widehat U_{x}.
\end{equation}

Let us continue the proof of the theorem. Consider the holomorphic function
$$
f_{y}:=\frac{\phi^{2}}{\phi^{2}(y)},\ \ \ y\in\widehat U_{x},
$$
defined in a neighbourhood of the closure of $\widetilde Y_{1}$. Here
$\phi$ is the same as in Lemma \ref{le3.3} and 
$\widehat U_{x}\subset\subset Y$ is the coordinate neighbourhood of $x\in K$
defined by (\ref{comp}). Then $f_{y}$
has double zeros at all $x_{j}$, $1\leq j\leq l$, is nonzero outside 
these points, and $f_{y}(y)=1$.

Finally, we introduce
\begin{equation}\label{e3.18}
F_{s,y}(w):=h_{s,y}^{2}(w)\cdot (r_{1}^{*}f_{y})(w)\ \ \ {\rm for\ all}\ \ \
w\in X,\ s\in S,\ y\in\widehat U_{x}.
\end{equation}
According to (\ref{e3.16''}) we have
\begin{equation}\label{e3.23}
\sum_{s\in S}|F_{s,y}(z)|\leq C'\ \ \ {\rm for\ all}\ \ \ z\in X,\ 
y\in\widehat U_{x},
\end{equation}
with $C':=\frac{C}{|\phi(y)|^{2}}$ depending on $x$, $Y$ only. 
Moreover,
\begin{equation}\label{e3.24}
F_{s,y}(y,t)=\delta_{st}
\end{equation}
where $\delta_{st}=0$ for $s\neq t$ and $\delta_{ss}=1$.

Using the functions $F_{s,y}$, $s\in S$, let us prove that
$r^{-1}(y)$ is interpolating for $H^{\infty}(X)$ for all $y\in\widehat U_{x}$.

In fact, for $a=\{a_{s}\}_{s\in S}\in l^{\infty}(S)$ and
$y\in\widehat U_{x}$ consider the function
\begin{equation}\label{e3.25}
[L_{y}(a)](z):=\sum_{s\in S}a_{s}F_{s,y}(z),\ \ \ z\in X.
\end{equation}
According to (\ref{e3.23}) we have
$$
\sup_{z\in X}|[L_{y}(a)](z)|\leq ||a||_{l^{\infty}}\cdot\sup_{z\in X}
\left(\sum_{s\in S}
|F_{s,y}(z)|\right)\leq C'||a||_{l^{\infty}}.
$$
Thus $L_{y}$ is a linear continuous 
operator from $l^{\infty}(r^{-1}(y))$ to $H^{\infty}(X)$ depending 
holomorphically on $y\in\widehat U_{x}$ with the norm bounded by a number
depending on $x$ and $Y$ only. Also, from (\ref{e3.24}) we obtain
\begin{equation}\label{e3.26}
[L_{y}(a)](y,t):=\sum_{s\in S}a_{s}F_{s,y}(y,t)=\sum_{s\in S}a_{s}
\delta_{st}=a_{t}=:a(y,t).
\end{equation}
That is, $L_{y}(a)|_{r^{-1}(y)}=a$. Therefore 
$r^{-1}(y)$, $y\in\widehat U_{x}$, is an interpolating sequence for
$H^{\infty}(X)$ with the constant of interpolation depending on $Y$ and $x$
only. Taking a finite open cover
of $K$ by sets $\widehat U_{x}$ and considering the corresponding interpolation
operators $L_{y}$ on $\widehat U_{x}$
we obtain that the constant of interpolation of each
$r^{-1}(x)$, $x\in K$, is bounded by a number depending on $K$ and $Y$ only.

The proof of Theorem \ref{te3} is complete.\ \ \ \ \ $\Box$
\sect{\hspace*{-1em}. Proof of Theorem \ref{te4}.}
\quad 
Let $\eta$ be a smooth bounded $(0,1)$-form on $X$ with 
$r(supp\ \eta)\subset K$ for some compact $K\subset\subset Y$. We must find
a smooth function $f$ on $X$ such that
\begin{equation}\label{e4.1}
\overline\partial f=\eta\ \ \ {\rm and}\ \ \
||f||_{L^{\infty}}:=\sup_{z\in X}|f(z)|\leq C||\eta||
\end{equation}
with $C$ depending on $K$, $Y$ and a hermitian metric $h_{Y}$ used
in the definition of the norm of $\eta$, see (\ref{e3}).

Without loss of generality we may assume that $\eta$ has compact
support. Indeed, let $\{X_{i}\}_{i\in\N}$, $X_{i}\subset\subset X$, be
an exhaustion of $X$ by relatively compact open domains. Let
$\{\chi_{i}\}_{i\in\N}$ be a family of smooth functions on $X$ such that
$\chi_{i}$ equals 1 on a subdomain $Z_{i}\subset\subset X_{i}$,
$0$ outside $X_{i}$ and $0\leq\chi_{i}\leq 1$, $i\in\N$. Assume also that
$\{Z_{i}\}_{i\in\N}$ forms an exhaustion of $X$, as well. Now, we set
$\eta_{i}:=\chi_{i}\cdot\eta$. 
Then $\{\eta_{i}\}$ converges to $\eta$ uniformly on compact
subsets of $X$ and $||\eta_{i}||\leq ||\eta||$ for all $i$. If we will
find smooth functions $f_{i}$ on $X$ satisfying the corresponding
conditions (\ref{e4.1}), then a standard normal family argument will give us
a subsequence $\{f_{i_{k}}\}_{k\in\N}$ of $\{f_{i}\}$ converging uniformly on
compact subsets of $X$ to a smooth function $f$ satisfying (\ref{e4.1}).
Thus it suffices to prove the theorem for the forms $\eta$ with compact 
supports.

Next, consider a finite open cover $(U_{i})_{1\leq i\leq n}$ of 
$K\subset\subset Y$ by sets $U_{i}:=\widehat U_{x_{i}}$, $x_{i}\in K$, defined
by (\ref{comp}).
By definition $(U_{i})_{1\leq i\leq n}$ 
also covers a neighbourhood $N\subset\subset Y$ of
$K$. Now we consider a finite open cover $(U_{i})_{n+1\leq i\leq m}$ of
$R\setminus N$ (where $R$ is a compact Riemann surface 
from (\ref{eq2.1}) containing $Y$) by coordinate
disks $U_{i}$ such that $U_{i}\cap K=\emptyset$ for
all $n+1\leq i\leq m$. Let $\{\rho_{i}\}_{1\leq i\leq m}$ be a smooth
partition of unity subordinate to the cover $(U_{i})_{1\leq i\leq m}$ of $R$.
Then since
$r(supp\ \eta)\subset K$, $U_{i}\cap K=\emptyset$ for all $n+1\leq i\leq m$,
and $supp\ \rho_{i}\subset U_{i}$ for $1\leq i\leq n$,
$$
\eta=\sum_{i=1}^{m}(r_{1}^{*}\rho_{i})\eta=\sum_{i=1}^{n}
(r_{1}^{*}\rho_{i})\eta.
$$
By the definition each $\eta_{i}:=(r_{1}^{*}\rho_{i})\eta$ is a smooth 
$(0,1)$-form with compact support
such that $r(supp\ \eta_{i})\subset U_{i}$. 
It suffices to prove the theorem for such forms $\eta_{i}$, i.e., 
to find smooth functions $f_{i}$ such that $\overline\partial f_{i}=\eta_{i}$
and $||f_{i}||_{L^{\infty}}\leq C_{i}||\eta_{i}|| \ (\leq C_{i}||\eta||)$
with $C_{i}$ depending on $U_{i}$, $Y$ and $h_{Y}$. Then 
$f:=\sum_{i=1}^{n}f_{i}$ satisfies the required statement of the theorem.

Thus without loss of generality we may assume that 
$r(supp\ \eta)\subset\subset \widehat U_{x}:=U$ for some $x\in K$ and
$\eta$ has compact support.
As before we identify $r^{-1}(U)$ with $U\times S$ where $S$ is the
fibre of $r$. Then there is a finite subset $S_{\eta}\subset S$ such that
$supp\ \eta\subset\subset U\times S_{\eta}$.
In a complex coordinate $z$ on $U$ the form $\eta$ is written as
$$
\eta(z,s)=g(z,s)\ \!d\overline{z},\ \ \ (z,s)\in U\times S_{\eta},
\ \ \ {\rm and}\ \ \ \eta=0\ \ \ {\rm outside}\ \ \ U\times S_{\eta}.
$$
By the hypothesis of the theorem we have
\begin{equation}\label{e4.2}
|\eta|:=\sup_{z\in U, s\in S_{\eta}}|g(z,s)|\leq c||\eta||
\end{equation}
for some $c$ depending on $U$ and $h_{Y}$ only.

The remaining part of the 
proof repeats literally the proof of Proposition 5.1 of [Br4].
We refer to this paper for details.

Consider the family of interpolation operators
$L_{z}: l^{\infty}(r^{-1}(y))\to H^{\infty}(X)$ 
holomorphic in $z\in U$ with norms bounded by a number
$C'$ depending on $U$ and $Y$ only, see (\ref{e3.25}). Let us define
the $(0,1)$-form $\lambda$ on $R$ with values in $H^{\infty}(X)$ 
by the formula
\begin{equation}\label{e4.3}
\lambda(z):=L_{z}(g(z,\cdot))\ \!d\overline{z},\ \ \ z\in U,\ \ \ {\rm and}
\ \ \ \lambda=0\ \ \ {\rm outside}\ \ \ U.
\end{equation}
Since $supp\ \eta\subset\subset U\times S_{\eta}$ and $S_{\eta}$ is a 
finite subset of $S$, the definition of $L_{z}$ in (\ref{e3.26})
implies that $\lambda$ is smooth. Using an integral formula we can
solve the equation $\overline\partial F=\lambda$ on a fixed
neighborhood of the closure of $Y$ in $R$ to get a smooth 
solution $F: Y\to H^{\infty}(X)$ satisfying 
\begin{equation}\label{e4.4}
\sup_{z\in Y,w\in X}|[F(z)](w)|\leq c'|\lambda|:=c'\cdot\sup_{z\in Y, w\in X}
|[L_{z}(g(z,\cdot))](w)|.
\end{equation}
with $c'$ depending on $Y$ only. Finally, we set
\begin{equation}\label{e4.5}
f(w):=[F(r(w))](w),\ \ \ w\in X.
\end{equation}
Since $\{L_{z}\}$ are interpolation operators holomorphic in $z\in U$ one has
$$
\overline\partial f(w):=[\lambda(r(w))](w)=\eta(w),\ \ \ w\in X.
$$
Moreover, from estimates (\ref{e4.2}), (\ref{e4.4}) and
$||L_{z}||\leq C'$ we get
$$
\sup_{w\in X}|f(w)|\leq C||\eta||
$$
where $C:=c'\cdot C'\cdot c$.

This completes the proof of the theorem.\ \ \ \ \ $\Box$ 
\sect{\hspace*{-1em}. Proof of Theorem \ref{te1}.}
\quad Let $r:X\to Y$ be an unbranched covering of a Caratheodory hyperbolic
Riemann surface of finite type $Y$. The fact that $X$ is Caratheodory 
hyperbolic follows easily from Theorem \ref{te3} and the 
Caratheodory hyperbolicity
of $Y$. Let us prove now the corona theorem for $H^{\infty}(X)$.

First we consider a finite open cover
${\cal U}=(U_{j})_{0\leq j\leq l}$ of 
$\widetilde Y:=Y\cup (\cup_{1\leq j\leq l}\{x_{j}\})$ such that for 
$1\leq j\leq l$ the set
$U_{j}\subset\subset\widetilde Y$ is an open coordinate disk 
centered at $x_{j}$  and $U_{0}$ is a 
bordered Riemann surface which intersects each $U_{j}$, $1\leq j\leq l$, 
by a set biholomorphic to an open annulus, see definition (\ref{eq2.1}). 
By $\{\rho_{j}\}_{0\leq j\leq l}$ we denote a smooth partition of unity on
$\widetilde Y$ subordinate to the cover ${\cal U}$.
We set 
$U_{j}^{*}:=U_{j}\setminus\{x_{j}\}$, $1\leq j\leq l$. Then $U_{j}^{*}$ is
biholomorphic to a punctured open disk in $\Co$. Now, $r^{-1}(U_{j}^{*})$
is a disjoint union of sets biholomorphic to $\Di$ or to the punctured disk
$\Di^{*}$ (because the fundamental group of $U_{j}^{*}$ is $\Z$).
Moreover, according to our construction $\pi_{1}(U_{0})\cong\pi_{1}(Y)$.
Hence $r^{-1}(U_{0})$ is an open connected subset of $X$.

Suppose now that a collection $f_{1},\dots, f_{n}$ of functions from
$H^{\infty}(X)$ satisfies the corona condition (\ref{e1}).
Since each connected component of $r^{-1}(U_{j}^{*})$ is biholomorphic to
$\Di$ or $\Di^{*}$, according to the Carleson Corona Theorem, see, e.g., 
[Ga, Ch.$\!$ VIII, Th.$\!$ 2.1], there are a constant $C_{1}(n,\delta)$
(with $\delta$ from (\ref{e1})) and functions $g_{1}^{j},\dots, g_{n}^{j}$
from $H^{\infty}(r^{-1}(U_{j}^{*}))$, $1\leq j\leq l$, such that
\begin{equation}\label{e5.1}
\begin{array}{c}
\displaystyle
f_{1}g_{1}^{j}+\cdots +f_{n}g_{n}^{j}\equiv 1\ \ \ {\rm on}\ \ \
r^{-1}(U_{j}^{*})\ \ \ {\rm and}\\
\\
\displaystyle
||g_{k}^{j}||\leq C_{1}(n,\delta),\ \ \ 1\leq j\leq l,\ 1\leq k\leq n.
\end{array}
\end{equation} 
Also, since $U_{0}\subset\subset Y$ is a bordered Riemann surface, 
according to [Br1, Cor.$\!$ 1.6] there are a constant $C_{2}(Y,n,\delta)$ and
functions $g_{1}^{0},\dots, g_{n}^{0}$
from $H^{\infty}(r^{-1}(U_{0}))$ such that
\begin{equation}\label{e5.2}
\begin{array}{c}
\displaystyle
f_{1}g_{1}^{0}+\cdots +f_{n}g_{n}^{0}\equiv 1\ \ \ {\rm on}\ \ \
r^{-1}(U_{0})\ \ \ {\rm and}\\
\\
\displaystyle
||g_{k}^{0}||\leq C_{2}(Y,n,\delta),\ \ \ 1\leq k\leq n.
\end{array}
\end{equation} 
We set
\begin{equation}\label{e5.3}
h_{k}:=\sum_{j=0}^{l}(r^{*}\rho_{j})g_{k}^{j},\ \ \ 
1\leq k\leq n,\ 0\leq j\leq l.
\end{equation}
Since $supp\ \rho_{j}\subset\subset U_{j}$, $0\leq j\leq l$, $h_{k}$ are
smooth functions on $X$ such that
\begin{equation}\label{e5.4}
\begin{array}{c}
\displaystyle
f_{1}h_{1}+\cdots +f_{n}h_{n}\equiv 1\ \ \ {\rm on}\ \ \ X\ \ \ {\rm and}\\
\\
||h_{k}||_{L^{\infty}}\leq C_{3}(Y,n,\delta),\ \ \ 1\leq k\leq n.
\end{array}
\end{equation}

Next we will use a standard construction based on the Koszul complex,
see [Ga, Ch.$\!$\ VIII]. Namely, we write
\begin{equation}\label{e5.5}
\begin{array}{c}
\displaystyle
g_{j}(z)=h_{j}(z)+\sum_{k=1}^{n}a_{j,k}(z)f_{k}(z)\\
\\
\displaystyle
a_{j,k}(z)=b_{j,k}(z)-b_{k,j}(z)\ \ \ {\rm and}\\
\\
\displaystyle
\overline\partial b_{j,k}=h_{j}\cdot\overline\partial h_{k}=:\eta_{j,k},
\ \ \ j\neq k.
\end{array}
\end{equation}
According to (\ref{e5.3}) and (\ref{e5.4}) the smooth $(0,1)$-forms
$\eta_{j,k}$ on $X$ satisfy
$$
r(supp\ \eta_{j,k})\subset\subset U_{0}\ \ \ {\rm and}\ \ \
||\eta_{j,k}||\leq C_{4}(Y,n,\delta)\ \ \ {\rm for\ all}\ \ \ j,\ k.
$$
where $||\cdot||$ is defined with respect to a fixed hermitian metric
$h_{Y}$ on $Y$, see (\ref{e4}). Therefore by Theorem \ref{te4} there are
smooth functions $b_{j,k}$ on $X$ satisfying equations (\ref{e5.5}) such that
$$
||b_{j,k}||_{L^{\infty}}\leq C_{5}(Y,n,\delta)\ \ \ {\rm for\ all}\ \ \
j,\ k.
$$
Then the functions $g_{j}$ on $X$ belong to $H^{\infty}(X)$ and satisfy
\begin{equation}\label{e5.6}
\begin{array}{c}
\displaystyle
f_{1}g_{1}+\cdots +f_{n}g_{n}\equiv 1\ \ \ {\rm and}\\
\\
\displaystyle
||g_{j}||\leq C(Y,n,\delta)\ \ \ {\rm for\ all}\ \ \ j.
\end{array}
\end{equation}

This completes the proof of Theorem \ref{te1}.\ \ \ \ \ $\Box$

\end{document}